\documentclass{article}

\usepackage{graphicx,etoolbox}

\usepackage{hyperref}

\usepackage{hyperref}

\usepackage{amssymb}

\newcommand{\newc}{\newcommand}

\newc{\eqnoset}{\setcounter{equation}{0}}

\newcommand{\mref}[1]{(\ref{#1})}
\newcommand{\reflemm}[1]{Lemma~\ref{#1}}
\newcommand{\refrem}[1]{Remark~\ref{#1}}
\newcommand{\reftheo}[1]{Theorem~\ref{#1}}

\newcommand{\refcoro}[1]{Corollary~\ref{#1}}

\newcommand{\refsec}[1]{Section~\ref{#1}}

\newcommand{\beq}{\begin{equation}}
	\newcommand{\eeq}{\end{equation}}
\newcommand{\beqno}[1]{\begin{equation}\label{#1}}
	
	\newcommand{\barr}{\begin{array}}
		\newcommand{\earr}{\end{array}}
	
	\newc{\bearr}{\begin{eqnarray*}}
		\newc{\eearr}{\end{eqnarray*}}
	
	\newc{\bearrno}[1]{\begin{eqnarray}\label{#1}}
		\newc{\eearrno}{\end{eqnarray}}
	
	\newc{\non}{\nonumber}
	\newc{\nol}{\nonumber\nl}
	
	\newcommand{\bdes}{\begin{description}}
		\newcommand{\edes}{\end{description}}
	\newc{\benu}{\begin{enumerate}}
		\newc{\eenu}{\end{enumerate}}
	\newc{\btab}{\begin{tabular}}
		\newc{\etab}{\end{tabular}}

	\newtheorem{theorem}{Theorem}[section]
	\newtheorem{defi}[theorem]{Definition}
	\newtheorem{lemma}[theorem]{Lemma}
	\newtheorem{rem}[theorem]{Remark}
	\newtheorem{exam}[theorem]{Example}
	\newtheorem{propo}[theorem]{Proposition}
	\newtheorem{corol}[theorem]{Corollary}
	\newtheorem{conj}[theorem]{Conjecture}

	\newcommand{\btheo}[1]{\begin{theorem}\label{#1}}
		\newc{\brem}[1]{\begin{rem}\label{#1}\em}
			\newc{\bexam}[1]{\begin{exam}\label{#1}\em}
				\newc{\bdefi}[1]{\begin{defi}\label{#1}}
					\newcommand{\blemm}[1]{\begin{lemma}\label{#1}}
						\newcommand{\bprop}[1]{\begin{propo}\label{#1}}
							\newcommand{\bcoro}[1]{\begin{corol}\label{#1}}
								\newcommand{\btheoc}[1]{\begin{conj}\label{#1}}
									\newcommand{\etheo}{\end{theorem}}
								\newc{\etheoc}{\end{conj}}
							\newcommand{\elemm}{\end{lemma}}
						\newcommand{\eprop}{\end{propo}}
					\newcommand{\ecoro}{\end{corol}}
				\newc{\erem}{\end{rem}}
			\newc{\eexam}{\end{exam}}
		\newc{\edefi}{\end{defi}}
	
	\newc{\rmk}[1]{{\bf REMARK #1: }}
	\newc{\DN}[1]{{\bf DEFINITION #1: }}
	
	\newcommand{\bproof}{{\bf Proof:~~}}
	\newc{\eproof}{{\vrule height8pt width5pt depth0pt}\vspace{3mm}}
	
	\newc{\bfrac}[2]{\dspl{\frac{#1}{#2}}}

	\newc{\nid}{\noindent}
	
	\newcommand{\dspl}{\displaystyle}
	\newc{\grad}{\nabla}
	\newc{\Div}{\mbox{div}}
	\newc{\pdt}[1]{\dspl{\frac{\partial{#1}}{\partial t}}}
	\newc{\pdn}[1]{\dspl{\frac{\partial{#1}}{\partial \nu}}}
	\newc{\pdNi}[1]{\dspl{\frac{\partial{#1}}{\partial \mathcal{N}_i}}}
	\newc{\pD}[2]{\dspl{\frac{\partial{#1}}{\partial #2}}}
	\newc{\dt}{\dspl{\frac{d}{dt}}}
	\newc{\bdry}[1]{\mbox{$\partial #1$}}
	\newc{\sgn}{\mbox{sign}}
	
	\newc{\Hess}[1]{\frac{\partial^2 #1}{\pdh z_i \pdh z_j}}
	\newc{\hess}[1]{\partial^2 #1/\pdh z_i \pdh z_j}

	\newc{\ag}{\alpha}
	\newc{\bg}{\beta}
	\newc{\cg}{\gamma}\newc{\Cg}{\Gamma}
	\newc{\dg}{\delta}\newc{\Dg}{\Delta}
	\newc{\eg}{\varepsilon}
	\newc{\zg}{\zeta}
	\newc{\thg}{\theta}
	\newc{\llg}{\lambda}\newc{\LLg}{\Lambda}
	\newc{\kg}{\kappa}
	\newc{\rg}{\rho}
	\newc{\sg}{\sigma}\newc{\Sg}{\Sigma}
	\newc{\tg}{\tau}
	\newc{\fg}{\phi}\newc{\Fg}{\Phi}
	\newc{\vfg}{\varphi}
	\newc{\og}{\omega}\newc{\Og}{\Omega}
	\newc{\pdh}{\partial}
	
	\newc{\ccG}{{\cal G}}

	\newc{\ii}[1]{\int_{#1}}
	\newc{\iidx}[2]{{\dspl\int_{#1}~#2~dx}}
	\newc{\bii}[1]{{\dspl \ii{#1} }}
	\newc{\biii}[2]{{\dspl \iii{#1}{#2} }}
	\newc{\su}[2]{\sum_{#1}^{#2}}
	\newc{\bsu}[2]{{\dspl \su{#1}{#2} }}

	\newc{\biiom}[1]{{\dspl\int_{\bdrom}~ #1 ~d\sg}}
	\newc{\io}[1]{{\dspl\int_{\Og}~ #1 ~dx}}
	\newc{\bio}[1]{{\dspl\int_{\bdrom}~ #1 ~d\sg}}
	\newc{\bsir}{\bsu{i=1}{r}}
	\newc{\bsim}{\bsu{i=1}{m}}
	
	\newc{\iibr}[2]{\iidx{\bprw{#1}}{#2}}
	\newc{\Intbr}[1]{\iibr{R}{#1}}
	\newc{\intbr}[1]{\iibr{\rg}{#1}}
	\newc{\intt}[3]{\int_{#1}^{#2}\int_\Og~#3~dxdt}
	
	\newc{\itQ}[2]{\dspl{\int\hspace{-2.5mm}\int_{#1}~#2~dz}}
	\newc{\mitQ}[2]{\dspl{\rule[1mm]{4mm}{.3mm}\hspace{-5.3mm}\int\hspace{-2.5mm}\int_{#1}~#2~dz}}
	\newc{\mitQQ}[3]{\dspl{\rule[1mm]{4mm}{.3mm}\hspace{-5.3mm}\int\hspace{-2.5mm}\int_{#1}~#2~#3}}
	
	\newc{\mitx}[2]{\dspl{\rule[1mm]{3mm}{.3mm}\hspace{-4mm}\int_{#1}~#2~dx}}
	\newc{\mitmu}[2]{\dspl{\rule[1mm]{3mm}{.3mm}\hspace{-4mm}\int_{#1}~#2~d\mu}}
	\newc{\iidmu}[2]{\iidx{#1}{#2}}
	
	\newc{\iidm}[3]{{\dspl\int_{#1}~#2~d #3}}
	
	\newc{\itQmu}[2]{\dspl{\int\hspace{-2.5mm}\int_{#1}~#2~d\mu}}
	\newc{\mitQmu}[2]{\dspl{\rule[1mm]{4mm}{.3mm}\hspace{-5.3mm}\int\hspace{-2.5mm}\int_{#1}~#2~d\mu}}

	\newc{\mitQq}[2]{\dspl{\rule[1mm]{4mm}{.3mm}\hspace{-5.3mm}\int\hspace{-2.5mm}\int_{#1}~#2~d\bar{z}}}
	\newc{\itQq}[2]{\dspl{\int\hspace{-2.5mm}\int_{#1}~#2~d\bar{z}}}

	\newc{\pder}[2]{\dspl{\frac{\partial #1}{\partial #2}}}

	\newc{\bdrom}{\bdry{\Og}}

	\newc{\bilhom}{\mbox{Bil}(\mbox{Hom}(\RR^{nm},\RR^{nm}))}
	\newc{\VV}[1]{{V(Q_{#1})}}
	
	\newc{\ccA}{{\mathcal A}}
	\newc{\ccB}{{\mathcal B}}
	\newc{\ccC}{{\mathcal C}}
	\newc{\ccD}{{\mathcal D}}
	\newc{\ccE}{{\mathcal E}}
	\newc{\ccH}{\mathcal{H}}
	\newc{\ccF}{\mathcal{F}}
	\newc{\ccI}{{\mathcal I}}
	\newc{\ccJ}{{\mathcal J}}
	\newc{\ccK}{{\mathcal K}}
	\newc{\ccP}{{\mathcal P}}
	\newc{\ccQ}{{\mathcal Q}}
	\newc{\ccR}{{\mathcal R}}
	\newc{\ccS}{{\mathcal S}}
	\newc{\ccT}{{\mathcal T}}
	\newc{\ccX}{{\mathcal X}}
	\newc{\ccY}{{\mathcal Y}}
	\newc{\ccZ}{{\mathcal Z}}
	
	\newc{\bb}[1]{{\mathbf #1}}
	
	\newc{\myprod}[1]{\langle #1 \rangle}
	\newc{\mypar}[1]{\left( #1 \right)}
	
	\newc{\BLLg}{\mathbf{\LLg}}
	
	\newc{\mA}{\mathbf{A}}
	\newc{\mB}{\mathbf{B}}
	\newc{\mC}{\mathbf{C}}
	\newc{\mD}{\mathbf{D}}
	\newc{\mE}{\mathbf{E}}
	\newc{\mF}{\mathbf{F}}
	\newc{\mJ}{\mathbf{J}}
	\newc{\mG}{\mathbf{G}}
	\newc{\mP}{\mathbf{P}}
	\newc{\mR}{\mathbf{R}}
	\newc{\mQ}{\mathbf{Q}}
	\newc{\mX}{\mathbf{X}}
	\newc{\muu}{\mathbf{u}}
	\newc{\mvv}{\mathbf{v}}
	
	\newc{\mllg}{\mathbb{\lambda}}
	\newc{\mLLg}{\mathbf{\LLg}}

	\newc{\lspn}[2]{\mbox{$\| #1\|_{\Lsp{#2}}$}}
	\newc{\Lpn}[2]{\mbox{$\| #1\|_{#2}$}}
	\newc{\Hn}[1]{\mbox{$\| #1\|_{H^1(\Og)}$}}

	\newc{\mynorm}[2]{\| #1\|_{#2}}

	\newcommand{\RR}{{\rm I\kern -1.6pt{\rm R}}}

	\newc{\itQQ}[2]{\dspl{\int_{#1}#2\,dz}}
	\newc{\mmitQQ}[2]{\dspl{\rule[1mm]{4mm}{.3mm}\hspace{-4.3mm}\int_{#1}~#2~dz}}
	\newc{\MmitQQ}[2]{\dspl{\rule[1mm]{4mm}{.3mm}\hspace{-4.3mm}\int_{#1}~#2~d\mu}}

	\newc{\MUmitQQ}[3]{\dspl{\rule[1mm]{4mm}{.3mm}\hspace{-4.3mm}\int_{#1}~#2~d#3}}
	\newc{\MUitQQ}[3]{\dspl{\int_{#1}~#2~d#3}}

	\newc{\mccP}{\mathbb{P}}
	\newc{\mccK}{\mathbb{K}}
	
	\newc{\DKTmU}{\mccK(U)}
	\newc{\DKTmUold}{(K_U(U)^{-1})^T}
	
	\newc{\myPi}{\mathbf{W}}
	\newc{\myIbar}{\bar{\ccI}_1}
	\newc{\myIhat}{\hat{\ccI}_1}
	\newc{\myIbreve}{\breve{\ccI}_0}
	
	\newc{\mmk}{\mathbf{k}}

	\newcommand{\ma}{\mathbf{a}}
	\newcommand{\mg}{\mathbf{g}}

	\newc{\mfu}{\mathbf{f_u}}
	\newc{\mh}{\mathbf{h}}
	\newc{\mb}{\mathbf{b}}
	\newc{\mf}{\mathbf{f}}

	\newcommand{\barrl}[2]{\barr{ll}\lefteqn{#1}\hspace{#2}&\\}

	\newc{\twomatrix}[1]{\left[\barr{cc}#1\earr\right]}
	\newc{\threematrix}[1]{\left[\barr{ccc}#1\earr\right]}

	\newc{\mN}{\mathbf{N}}
	\newc{\mI}{\mathbf{I}}
	\newc{\mH}{\mathbf{H}}

	\newc{\mk}{\mathbf{k}}
	\newc{\mr}{\mathbf{r}}

	\newc{\DIAGM}[2]{\left[\barr{ccc}#1&0\ldots&0\\
		\vdots&\ddots&\vdots\\0&\ldots0&#2\earr \right]}
	\newc{\DiagM}[2]{\mbox{diag}\left[#1
		\cdots #2 \right]}
	\newc{\vVEC}[2]{\left[\barr{c}#1\\
		\vdots\\#2\earr \right]}
	\newc{\hVEC}[2]{\left[#1
		\cdots #2 \right]}
	
	\newc{\mq}{\mathbf{q}}

	\newc{\msys}[1]{\left\{\barr{l}#1\earr
		\right.}
	\newc{\msysa}[1]{\left\{\barr{ll}#1\earr
		\right.}
	
	\newc{\bbM}{\mathbb{M}}
	\newc{\mat}[1]{\left[\barr{cc}#1\earr\right]}
	
	\newc{\me}{\mathbf{e}}

	\newc{\vecc}[2]{\left[\barr{cc}#1\\#2\earr\right]}
	\newc{\mL}{\mathbb{L}}
	
	\newc{\cO}{{\cal O}}
	
	\newc{\cM}{{\cal M}}
	
	\newc{\myega }{\eg_0(R)}
	\newc{\myeg}{\eg_1(\eg_*)}
	\newc{\myegp}{\hat{\eg}_1(\eg_*)}

\newc{\diagA}{\mathbb{A}_d}
\newc{\mBB}{\mathbb{B}}
\newc{\MLT}[1]{{\cal M}_{lt}(\Og,#1)}
\newc{\ALT}[1]{{\cal A}_{l}(\Og,#1)}
\newc{\mM}{\mathbb{M}}

\newc{\diag}[1]{\mbox{diag}(#1)}
\newc{\off}[1]{\mbox{offdiag}(#1)}
\newc{\mT}{\mathbb{T}}

\usepackage{amssymb}

\begin{document}

	\vspace*{-.8in}
	\begin{center} {\LARGE\em Regularity of Solutions to a Class of Degenerate Cross Diffusion Systems of $m$ Equations .}
		
	\end{center}

	\vspace{.1in}
	
	\begin{center}

		{\sc Dung Le}{\footnote {Department of Mathematics, University of
				Texas at San
				Antonio, One UTSA Circle, San Antonio, TX 78249. {\tt Email: Dung.Le@utsa.edu}
				}}

	\end{center}

	\begin{abstract} We study the regularity of weak solutions and global existence of classical  to cross diffusion systems of porous media type of $m$ equations on $N$-dimensional domains ($m,N\ge2$).
	 \end{abstract}

	\section{Introduction}\label{intro}
	
	This paper deals with the cross diffusion system
	\beqno{exsyszpara}W_t=\Div(A(W)DW)+G(W)W\eeq
	on bounded domain $\Og$ of $\RR^N$ with homogeneous Dirichlet or Neumann boundary conditions. Here, $W$ is an unknown vector in $\RR^m$, $m\ge 2$. As usual, $A(W), G(W)$ are $m\times m$ matrices and  for some positive functions $\llg(W),\LLg(W)$ we assume the ellipticity condition
	\beqno{normalellcond}\LLg(W)|DW|^2\ge \myprod{A(W)DW,DW}\ge \llg(W)|DW|^2.\eeq
	
	We will present several a priori estimates for higher norms strong and weak solutions of the above  system by weaker norms. These estimates play crucial roles in the regularity and  global existence problems of \mref{exsyszpara}.
	
	The regularity of weak solutions to (regular or non-regular) \mref{exsyszpara} is a long standing problem in the theory of partial differential systems (e.g. see \cite{GiaS}). It is well known that  weak solutions are H\"older continuous around points where the BMO norms are small in small balls.
	
	On the other hand, the global existence of strong/classical solutions has been considered in a pioneering paper of Amann \cite{Am2} where he showed that a solution would exist globally if its $W^{1,p}(\Og)$ norm for some $p>N$ does not blow up in finite time.
	
	Recently, in \cite{dlebook,dlebook1}, we introduce the strong/weak Gagliardo-Nirenberg  inequalities involving BMO norms (GNBMO for short) and assume that $A(W)$ satisfies  the spectral gap condition to prove that  the smallness of BMO norms in small balls also yields a control on the $W^{1,p}(\Og)$ norm of weak solutions and thus they are classical and exist globally.

	In \refsec{degGNsec} of this paper, we will show that new weighted versions of the strong/weak GNBMO inequalities are also available. We will sketch the proof of  strong/weak GNBMO inequalities and present some simple modifications to obtain new weighted versions. The latter was introduced to deal with degenerate systems as the regularity results of \cite{GiaS}  cannot be used for degenerate systems if we assume only that $\llg(W)>0$. In \refsec{regsyssec} we will study  these systems and show that they possess weak solutions that are H\"older continuous. 
	Again, we can see the involvement of BMO norms, the property that the BMO norms are small is critical in the proofs of regularity and global existence of solutions.
	The smallness of BMO can be established for thin domains.
	
	Uniqueness and boundedness of weak solutions will be considered in \refsec{unisec}. We will study a class of porous media type systems and improve some results in \cite{Vasquez}.

	
	
	\section{Weighted Gagliardo-Nirenberg inequalities involving BMO norms}\label{degGNsec}\eqnoset

	In this section we revisit the key ingredients of the proof of our existence theory and regularity properties of solutions to non-degenerate systems in  \cite{dlebook1}: the local weighted Gagliardo-Nirenberg inequalities involving BMO norm. We will establish new versions of these inequalities which can apply to degenerate systems. 
	
	Firstly, let us recall some notations.
	For any measurable subset $A$ of $\Og$  and any  locally integrable function $U:\Og\to\RR^m$ we denote by  $|A|$ the measure of $A$ and $U_A$ the average of $U$ over $A$. That is, $$U_A=\mitx{A}{U(x)} =\frac{1}{|A|}\iidx{A}{U(x)}.$$
	
	We also recall some well known notions from Harmonic Analysis.

	A function $f\in L^1(\Og)$ is said to be in $BMO(\Og)$ if \beqno{bmodef} [f]_{*}:=\sup_Q\mitx{Q}{|f-f_Q|}<\infty.\eeq We then define $$\|f\|_{BMO(\Og)}:=[f]_{*}+\|f\|_{L^1(\Og)}.$$

	For $\cg\in(1,\infty)$ we say that a nonnegative locally integrable function $w$ belongs to the class $A_\cg$ or $w$ is an $A_\cg$ weight if the quantity
	\beqno{aweight} [w]_{\cg} := \sup_{B\subset\Og} \left(\mitx{B}{w}\right) \left(\mitx{B}{w^{1-\cg'}}\right)^{\cg-1} \quad\mbox{is finite}.\eeq
	Here, $\cg'=\cg/(\cg-1)$ and the supremum is taken over all cubes $B$ in $\Og$. For more details on these classes we refer the reader to \cite{FPW,OP}.  We denote by $l(B)$ the side length of $B_{l(B)}$ and by $ B_{\tau l(B)}$ the cube which is concentric with $B$ and has side length $\tau l(B)$. Also, we write $\Og_R=\Og\cap B_R$.
	
	We will only sketch the main modifications needed to prove the new inequality with weights to cover the degenerate cases. We refer the readers to \cite[section 2.4]{dlebook1} for more details.

	\subsection{The sketch of the proof for the {\em weighted} strong Gagliardo-Nirenberg inequality with BMO norms:} \label{GNineqsec}\index{The strong Gagliardo-Nirenberg inequality}

	We will always use the following notations and hypotheses.
	Fixing $R>0$, we let $\og$ be a smooth cutoff function in $B_R$ ($|D\og|\le C/R$).
	 Let $\Gamma, \llg$ be local bounded nonnegative functions on $\Og$. We consider the following integrals (compared with $I_1,I_2$ in \cite{dlebook1}).
	
	\beqno{Idef} I_1:=\iidx{\Og}{\Gamma |Du|^{2p+2}\og^2},\;
	I_2:=\iidx{\Og}{\llg|Du|^{2p-2}|D^2u|^2\og^2},\eeq
	\beqno{Idef2}\myIbreve:=\iidx{\Og}{\llg|Du|^{2p}\og^2}.\eeq

	We modify the proof of \cite[Theorem 2.4.1]{dlebook1} to establish the following inequality, which will be referred to as the weighted strong GNBMO inequality in the sequel. This is the {\em strong} version because of the involvement of the second derivative $D^2u$ and we will present its {\em weak} version later.
	
	\btheo{GNlobal}   
	Suppose that  $\myega,I_1,I_2,\myIbreve$ are finite and that there is a constant $C_*(R)$ such that
	\beqno{Gammallg} \Gamma\le C_*\llg \mbox{ and } |D\Gamma|^2|u|^2\le C\Gamma\llg.\eeq
	
	 Define $\myega =\|u\|_{BMO(\Og_{2R})}$. Then,   there are  constants $C(N), C(N,p)$ such that \beqno{wGNglobalest}I_1\le C(N)C_*^2\myega ^2 I_2+[C(N,p)+\frac{C_*^2\myega ^2}{R^2}]\myIbreve.\eeq
	\etheo

	We note the difference here and the old versions in \cite{dlebook1}: the condition \mref{Gammallg} and the assertion \mref{wGNglobalest} which involve with $\Gamma,\llg,C_*,C(N,p)$. This is needed for our purpose later.

	Again, one of the key ingredients of our proof is the duality between the Hardy space $\ccH^1(\Og)$ and $BMO(\Og)$ space. This is the famous Fefferman-Stein theorem (see \cite{st}). It is well known that the norm of the Hardy space can be defined by
	$$\|F\|_{\ccH^1(\Og)}=\|F\|_{L^1(\Og)}+\|M_*F\|_{L^1(\Og)},\; M_{*}F(y) =\sup_{\fg}\int_{\Og}\fg(x-y)f(x)dx.$$ Here, the supremum is taken over all $\fg\in C^1(\RR^N)$ which has support in some cube $Q\subset \RR^N$ centered at $y$ with side length $l(Q)$ and satisfies $$0\le \fg(x)\le l(Q)^{-N},\; |D\fg(x)|\le l(Q)^{-N-1} \mbox{ for all $x\in Q$}.$$

	We will also use the definition of the  {\em centered}  Hardy-Littlewood maximal operator acting on functions $F\in L^1_{loc}(\Og)$
	\beqno{maximal} M(F)(y) = \sup_\eg\{\mitx{B_\eg(y)}{F(x)}\,:\, \eg>0 \mbox{ and } B_\eg(y)\subset\Og\}.\eeq
	
	We also  apply the Muckenhoupt theorems for non doubling measures. By \cite[Theorem 3.1]{OP}, we have that if $w$ is an $A_q$ weight then for any $F\in L^q(\Og)$ with $q>1$
	\beqno{HL1}\iidx{\Og}{M(F)^q w} \le C([w]_q)\iidx{\Og}{F^q w}.\eeq
	In particular, \beqno{HL}\iidx{\Og}{M(F)^q} \le C\iidx{\Og}{F^q}.\eeq

	\bproof (Sketch of proof for \reftheo{GNlobal})
	We start with the proof of \cite[Theorem 2.4.1]{dlebook1}. Denote $H=Du$, we can assume that $H$ is continuous as the inequality can be obtained by approximation. Applying
	integration by parts and noting that the boundary integral is zero, we then have 
	$$\barr{lll}I_1&=&-\iidx{\Og}{\Div(\Gamma |H|^{2p}Du)u\og^2}-\iidx{\Og}{2\Gamma\myprod{|H|^{2p}Du,D\og}u\og}\\
	&=& -\iidx{\Og}{\myprod{D\Gamma, |H|^{2p}Du)u\og^2}}-\iidx{\Og}{\Gamma \Div(|H|^{2p}Du)u\og^2}-\iidx{\Og}{2\Gamma\myprod{|H|^{2p}Du,D\og}u\og}.\earr$$
	
	So, with $G_1,G_2$ are similar to those in the proof of \cite[Theorem 2.4.1]{dlebook1} with the extra factor $\Gamma$ in their integrands,
	we can write
	$$I_1=-G_0-G_1-G_2,\; G_0=\iidx{\Og}{\myprod{D\Gamma, |H|^{2p}Du)u\og^2}}.$$ 
	
	We first consider $G_1$ and show that  $g=\Gamma\Div(|H|^{2p}Du)$ belongs to the Hardy space ${\cal H}^{1}$ and satisfies
	\beqno{gh1est}\|g\|_{{\cal H}^1} = \iidx{\Og}{\sup_\eg|g*\fg_\eg|} \le C_* I_2^\frac12I_1^\frac12+\frac{C_*}{R}I_1^\frac{1}{2}\myIbreve^\frac12.\eeq
	Once this is established, the Fefferman-Stein theorem on the duality of the BMO and Hardy spaces yields  $|G_1|\le \myega \|g\|_{{\cal H}^1}$ ($\myega =\|u\|_{BMO(\Og_R)}$) so that together with the estimate for $G_2$ (which will be estimated similarly) we will have
	\beqno{keyGN}I_1 \le \myega \|g\|_{{\cal H}^1}+\frac{C_*}{R}I_1^\frac{1}{2}\myIbreve^\frac12\le C_*\myega I_1^\frac{1}{2}I_2^\frac{1}{2}+\frac{C_*}{R}I_1^\frac{1}{2}\myIbreve^\frac12-G_0.\eeq

	As in \cite{dlebook1}, we write $g=g_1+g_2$ with $g_i=\Gamma\Div (V_i)\og^2$, setting  $h=|H|^{p-1}H$ and
	$$V_1= hDu\left(h-\mitx{B_\eg}{h}\right),\;V_2= hDu\mitx{B_\eg}{h}.$$
	
	Let us consider $g_1$ first. For any $y\in\Og$ and $B_\eg=B_\eg(y)\subset \Og$, we use integration by parts, ignoring $\og^2$ for simplicity  because the results involving $\og D\og$ can be treated the same way for $G_2$, and the property of $\fg_\eg$ and we can  exactly follow \cite{dlebook1} to have the following
	\beqno{hest}\frac{C_1}{\eg}
	\left(\mitx{B_\eg}{|h-h_{B_\eg}|^s}\right)^\frac1s \le C\Psi_2,\eeq where $\Psi_2(y)=\left(M(|H|^{(p-1)s_*}|DH|^{s_*})(y)\right)^\frac1{s_*}$.
	Setting $\Psi_3(y)=\left(M(|hDu|^{s_*})(y)\right)^\frac1{s_*}$ and putting these estimates together we thus have (see \cite[(2.31)]{dlebook1}) if \mref{Gammallg} holds  then $\Gamma\le C_*\llg$ so that
	\beqno{g1a}\sup_{\eg>0,}\Gamma|g_1*\fg_\eg| \le C(N)\Gamma
	\Psi_2\Psi_3 =C(N)\llg^\frac{1}{2}\Psi_2\frac{\Gamma}{\llg^\frac{1}{2}}\Psi_3\le C_*^\frac{1}{2}C(N)\llg^\frac{1}{2}\Psi_2\Gamma^\frac{1}{2}\Psi_3.\eeq

	Take $s=2n/(n-1)$ then $s_*=s'=2n/(n+1)$. With these notations  we can use \mref{HL1} with the measure $wdx=\og^2dx$,  because $2>2n/(n+1)=s_*$, to get
	$$\left(\iidx{\Og}{\llg\Psi_2^{2}\og^2}\right)^\frac1{2} \le \left\|M\left(\llg^\frac{s_*}{2}|H|^{(p-1)s_*}|DH|^{s_*}\right)\right\|_{L^\frac{2}{s_*}(\Og,\og^2dx)}^\frac{1}{s_*}\le C\|\llg^\frac{1}{2}|H|^{p-1}|DH|\|_{L^2(\Og,\og^2dx)}=I_2.$$
	
	Similarly, from the definitions of $I_1, h$, $$\left(\iidx{\Og}{\Gamma\Psi_3^2\og^2}\right)^\frac12\le C\|\Gamma^\frac{1}{2}|H|^{p}Du\|_{L^2(\Og,\og^2dx)} = CI_1^\frac12.$$

	Therefore, by Holder's inequality, the above estimates and the notations \mref{Idef} and \mref{Idef2} imply that \beqno{g1est}\iidx{\Og}{\sup_\eg|g_1*\fg_\eg|} \le C(N)C_*^\frac{1}{2} I_1^\frac12 I_2^\frac12.\eeq
	
	An estimate for $G_2$ can be obtained similarly.
	We consider $G_0$
	$$G_0=\iidx{\Og}{\myprod{D\Gamma, |H|^{2p}Du u\og^2}}\le \left(\iidx{\Og}{\frac{|D\Gamma|^2|u|^2}{\Gamma}|H|^{2p}\og^2}\right)^\frac{1}{2}\left(\iidx{\Og}{\Gamma|H|^{2p}|Du|^2\og^2}\right)^\frac{1}{2}.$$ 
	
	By \mref{Gammallg}, $|D\Gamma|^2|u|^2\le C\Gamma\llg$ so that
	$$G_0\le \left(\iidx{\Og}{\llg|H|^{2p}\og^2}\right)^\frac{1}{2}I_1^\frac{1}{2}= \myIbreve^\frac{1}{2}I_1^\frac{1}{2}.$$

	A simple use of Young's inequality to the right hand side then gives \mref{wGNglobalest}. \eproof
	
	\brem{Gammallgrem1} The theorem applies when $\Gamma,\llg$ are functions in $x$ in general so that we assume both conditions in \mref{Gammallg} As $\Gamma(u(x))\sim |D\Gamma(u(x))||u(x)|$ so that the first condition in \mref{Gammallg} is sufficient. Indeed, the first condition in \mref{Gammallg} implies the second and we have $|D\Gamma|^2|u|^2\le C\Gamma\llg$. It is easy to have $\myIbreve<\infty$ when $p=1$ so that we can start an induction argument  for bigger $p$. 	
	
	\erem

	\brem{Gammallgrem} If $\Gamma=\Gamma(u)$, a function in $u$ and $\myprod{\Gamma_u,u}\le C\Gamma(u)$ or $|\Gamma_u||u|\le CG(u)
	$ then 
	$$G_0=\iidx{\Og}{\myprod{D\Gamma, |H|^{2p}Duu\og^2}}=\iidx{\Og}{|H|^{2p}|Du|^2\myprod{\Gamma_u,u}\og^2}\le CI_1.$$
	Now, as $\myega =\|u\|_{BMO(\Og_{2R})}$, instead of \mref{wGNglobalest} we have $$I_1\le C(N)C_*^2\myega ^2 I_2+\frac{C_*^2\myega ^2}{R^2}\myIbreve +C\|u\|_{BMO(\Og_{2R})}I_1.$$
	If $\|u\|_{BMO(\Og_{2R})}<1/C$ then we have a usual estimate for $I_1$ by \mref{wGNglobalest}.
	We then need only the first condition in \mref{Gammallg}. We can apply this for porous media type systems and use approximations.
	
	\erem

	\brem{degsysrem} Usually, $\Gamma=\frac{|\ma_W|^2}{\llg}$ and $\llg$ is the smallest eigenvalue of $\ma(W)$. If $C_*\sim \|u\|_{BMO(\Og_{2R})}^{-2}$ then the first condition in \mref{Gammallg} is $\|u\|_{BMO(\Og_{2R})}^2|\ma_W|^2\le \llg^2$. This weaken a crucial condition in \cite{dlebook1} where we assumed that $\frac{|\ma_W|^2}{\llg}$ is bounded.
	
	\erem

	\subsection{A new (weak) {\em weighted} Gagliardo-Nirenberg inequality} \index{(weak) Gagliardo-Nirenberg inequality}We consider the case $H\ne Du$ and assume that $DH,Du$ are only sufficiently integrable. This is a quite improvement of \reftheo{GNlobal}. The new inequality presented here allows us to consider weak solutions and study their H\"older regularity. It resembles \mref{wGNglobalest} and the proof is also quite similar but needs some modifications as in the proof of \reftheo{GNlobal}. We only state the result here.
	
	Again, we can assume that $H,Du$ are continuous as the inequality can be obtained by approximation. We denote (here, $DH$ should be understood in its distribution sense).
	
	\newc{\mmyIbreve}{\mathbf{\myIbreve}}
	
	\beqno{Idefz} \mI_1:=\iidx{\Og}{\Gamma|H|^{2p}|Du|^2\og^2},\;
	\mI_2:=\iidx{\Og}{\llg|H|^{2p-2}|DH|^2\og^2},\eeq
	\beqno{Idef2z}\mmyIbreve:=\iidx{\Og}{\llg|H|^{2p}\og^2}.\eeq

	Especially, for a given $\eg_*>0$ we define
	\beqno{eg1def} \myeg=\left(\eg_*^{-N+2}\iidx{B_{\eg_*}}{|Du|^2}\right)^\frac12.\eeq

	We present here another version of \reftheo{GNlobal} which will be referred to as the weak GNBMO weighted inequality for short throughout this paper. Note also that we also have \refrem{Gammallgrem} and those follow it.
	
	\btheo{wGNlobalz}   Assume $\eg_*\le R$.
	Suppose that  $\myega, \mI_1,\mI_2,\mmyIbreve$ are finite and \mref{Gammallg}.
	Then   there is a  constant $C$ such that \beqno{GNglobalestz}\mI_1\le C(N)C_*^2\myega ^2 \mI_2+C(\frac{C_*^2\myega ^2}{R^2}+\frac{\myeg^2}{\eg_*^2}+1)\mmyIbreve.\eeq
	\etheo

	\section{Regularity for degenerate systems:}\label{regsyssec}\eqnoset The regularity results of \cite{GiaS}  cannot be used for degenerate systems (we now have $\llg(W)>0$ only).
	
	However, assuming the {\em spectral-gap condition} and arguing as in \cite{dlebook1}, we see that under the assumption that the BMO norms of a weak (weak-strong) solution $W$ in small balls in $\Og$ are small then we can use the same induction argument in \cite[Lemma 4.4.1 or Theorem 4.4.6]{dlebook1} and the weighted strong/weak Gagliardo-Nirenberg BMO inequalities (\reftheo{GNlobal} or \reftheo{wGNlobalz}) in finite steps to have a bound for $\|DW\|_{L^{2p}(B_R)}$ for some $p>N/2$ even if the system is degenerate. Of course, as $p>N/2$ then $W$ is H\"older continuous.
	
	We will present new regularity results for degenerate systems
	\beqno{degsys} W_t=\Div(\ma(W)DW)+F(W), \eeq
	where there are $\llg(W),\LLg(W)>0$ such that $\llg(W)|\zeta|^2\le \myprod{\ma(W)\zeta,\zeta}\le \LLg(W)|\zeta|^2$. Let $\nu_*=\sup_{W\in\RR^m}\frac{\llg(W)}{\LLg(W)}$.
	
	In fact, we can regularize \mref{degsys} by considering
	\beqno{degsysregularized} W_t=\Div([\ma(W)+\eg Id]DW)+F(W), \eeq
	where $\eg>0$. The estimates for $DW_\eg$ of {\em strong} solutions $W_\eg$ below (via \reftheo{GNlobal} and \refrem{Gammallgrem}) are uniform in $\eg$ so that we can take the limit $\eg\to0$.
	
	\btheo{regdegthm} Let  $W$ be a {\em strong} solution to \mref{degsysregularized}. Suppose  the spectral-gap condition $\nu_*>1-2/N$.  
	
	Importantly, we assume that  $\Gamma:=\frac{|\ma_W(W)|^2}{\llg(W)}$ and $\llg(W)$ satisfy \mref{Gammallg} of \reftheo{GNlobal} (see also \refrem{Gammallgrem}).
	Assume furthermore that for sufficiently small $R>0$ we have ($C(N)$ is the constant in \reftheo{GNlobal}) \beqno{BMOWcond}C(N)C_*^2\|W\|_{BMO(\Og_{2R})}^2<1.\eeq
	
	If there is $p>N/2$ such that $F_W\in L^{2p+2}(\Og)$, then for any $t_0>0$ there is a finite constant $C(p,t_0,\|W\|_{L^1(\Og)})$, {\em independent of $\eg$ in \mref{degsysregularized}}, such that $$\sup_{t\in(0,t_0)}\|DW\|_{L^{2p}(\Og_R)}\le C(p,t_0,\|F_W(W)\|_{L^{2p+2}(\Og)},\|W\|_{L^1(\Og)}).$$ Moreover, $W$ is H\"older continuous, $W\in C^\ag(\Og)$ for some $\ag>0$.
	
	Note that the {\em spectral-gap} condition is void if $m=1$ and we can have the result for all $\ag\in(0,1)$.

	\etheo
	
	\bproof We follow the proof of \cite[Lemma 4.4.1]{dlebook1} to estimate $\|DW\|_{L^{2p}_{loc}(\Og)}$. Applying the difference operator  $\dg_h$ in $x$ to  the equation of $W$, we see that $W$ weakly solves (dropping the terms involving $\mB$ in the proof of \cite[Lemma 4.4.1]{dlebook1})
	$$(\dg_hW)_t=\Div(\ma D(\dg_h W)+\ma_W \myprod{\dg_h W,DW})+F_W  \dg_hW.$$
	
	Let $p\ge1$ and   $\fg,\eta$ be positive $C^1$ cutoff functions for the concentric balls $B_s, B_t$ and the time interval $I$. For any $0<s<t<2R_0$ we test this system with $|\dg_hW|^{2p-2}\dg_hW\fg^2\eta$ and use Young's inequality for the term $|\ma_W||\dg_hW|^{2p-1}|DW||D(\dg_hW)|$ and the spectral gap condition (with $X=\dg_h W$). We get    for some constant $c_0>0$, $\Psi=|F_W|$,  and $Q=\Og_t\times I$
	$$\barrl{\sup_{t\in(t_0,T)}\iidx{\Og}{|\dg_hW|^{2p}\fg^2}+c_0\itQ{Q}{ \llg|\dg_hW|^{2p-2}|D(\dg_hW)|^2\fg^2}\le  }{1cm}& C\iidx{\Og\times\{t_0\}}{|\dg_hW|^{2p}\fg^2}+C\itQ{Q}{\frac{|\ma_W|^2}{\llg}|\dg_hW|^{2p}|DW|^2\fg^2}+\\&C\itQ{Q}{[\Psi|\dg_hW|^{2p}\fg^2+|\ma_W|(|\dg_hW|^{2p}|DW|)\fg|D\fg|]}.\earr$$
	
	We make use of \reftheo{GNlobal}  to estimate the second term on the right hand side, absorb it to the left  and proceed in the same way as in \cite[Lemma 4.4.1]{dlebook1} and obtain the estimate for $\sup_I\|DW\|_{L^{2p}(\Og_R)}$. \eproof
	
	\brem{strongvsweakGN} Here, we chose to use the strong version \reftheo{GNlobal} and deal with strong solutions of \mref{degsysregularized}. By approximation, i.e. letting $\eg\to0$, we find {\em a weak} solution of \mref{degsys} which are H\"older continuous. We can use the strong version \reftheo{wGNlobalz} and deal directly with weak solutions of \mref{degsys} but there degeneracy is quite weak. To avoid this, we choose to work with approximation as described (the uniqueness of weak solutions is another issue, if $f$ is Lipschitz in $W$ and $\ma(W)$ is mononicity in some sense then we have uniqueness, see \cite{DB}).
	
	\erem
	
	\brem{spectralgaprem} In many cases, by induction we can prove that $\|F_W(W)\|_{L^{2p+2}(\Og)}$ is bounded in terms of $\|W\|_{L^{1}(\Og)}$.
	We also note that the above argument holds as long as $\nu_*>1-1/p$. \erem

	We can deal directly with a {\em weak} solution to \mref{degsys} and follow the proof of \cite[Theorem 4.4.6]{dlebook1} by using \reftheo{wGNlobalz}.
	
	Concerning \refrem{Gammallgrem}, as $\Gamma=\frac{|\ma_W|^2}{\llg}$ is a function of $W$, if $\llg(W)\sim |W|^k$ for some $k\ge0$ then $\Gamma(W)$ is also a power of $|W|$ so that $\myprod{\Gamma_W,W}\le C\Gamma(W)$. 
	
	\brem{existenceweakrem} The condition \mref{BMOWcond} $C(N)C_*^2\|W\|_{BMO(\Og_{2R})}^2<1$ combined with the first one $\Gamma\le C_*\llg$ in \mref{Gammallg} applied to \mref{degsysregularized}  is now $\|W\|_{BMO(\Og_{2R})}|\ma_W|\le c(N)[|W|^k+\eg]$ or (as $|\ma_W(W)|\sim |W|^{k-1}$) \beqno{kpowercond}\|W\|_{BMO(\Og_{2R})}\le c(k)[|W|+\eg|W|^{1-k}]\eeq for some constant $c(k)>0$ small. For general $k,W$, this condition is hard to be achieved.
	
	Of course, \mref{kpowercond} can be verified if $\|W\|_{BMO(\Og_{2R})}$ is small, $N\ge2$ ($N>2$ and the domain is thin) and $k>1$ because the right hand side is bounded from below by a  positive constant $c$ (independent of $\eg$ as $\eg\to0$). If $k=1$ and $N=2$, $c$ depends on $\eg$ but we can follow the proof of \cite[Theorem 4.4.6]{dlebook1} to show $\|W\|_{BMO(\Og_{2R})}$ is small if $R$ small (the growth rates of $F,\llg$ in \refcoro{2Ndomaindeg} below allow us to obtain \cite[Theorem 3.1]{unib} when $N=2$). Meanwhile, we still derive that $\sup_I\|DW\|_{L^{2p}(\Og_R)}$ is still bounded  (by the constant $C(p,t_0,\|F_W(W)\|_{L^{2p+2}(\Og)},\|W\|_{L^1(\Og)})$ in \reftheo{regdegthm}) independent of $\eg$. This argument an be used to show the regularity of weak solutions by approximation ($\eg\to0$) if $1/C_*\sim \frac{(\llg(W)+\eg)^2}{|\ma_W(W)|^2}$ {\em is bounded from below by a positive number independent} of $\eg>0$ when $|W|$ small and $k\ge0$.
	
	\erem
	
	Hence, \reftheo{GNlobal} can be used in the proof and we have
	\bcoro{regdegco} Let  $W$ be a {\em strong} solution to \mref{degsysregularized}. Suppose  the spectral-gap condition $\nu_*>1-2/N$.  
	
	Importantly, we assume that $\llg(W)\sim |W|^k$ for some $k\ge0$ and $\|W\|_{BMO(\Og_{2R})}$ is sufficiently small.
	If there is $p>N/2$ such that if $F_W\in L^{2p+2}(\Og)$, then for any $t_0>0$ there is a finite constant $C(p,t_0,\|W\|_{L^1(\Og)})$, {\em independent of $\eg$ in \mref{degsysregularized}}, such that $$\sup_{t\in(0,t_0)}\|DW\|_{L^{2p}(\Og_R)}\le C(p,t_0,\|F_W(W)\|_{L^{2p+2}(\Og)},\|W\|_{L^1(\Og)}).$$ Moreover, $W$ is H\"older continuous, $W\in C^\ag(\Og)$ for some $\ag>0$ with uniform bounded norm.

	\ecoro
	
	If $\ma(W)=|W|^k Id$ for some $k\ge1$ then $\nu_*=1$ and the above theorem together with \refrem{spectralgaprem}
	show that $W\in C^\ag(\Og)$ for all $\ag\in(0,1)$.
	
	\brem{generaldeg} We see that for a more general case, if we have $\llg(W)>0$ and we can also prove $\|W\|_{BMO(\Og_{2R})}|\ma_W|\le c(N)[\llg(W)+\eg]$, then by the same argument to have $\sup_I\|DW\|_{L^{2p}(\Og_R)}$ is still bounded  (by the constant $C(p,t_0,\|W\|_{L^1(\Og)})$ in \reftheo{regdegthm}) independent of $\eg$.
	\erem

	Thus, the proof that $\|W\|_{BMO(\Og_{2R})}$ is small (independent of $\eg$) is important. For $N=2$, we can follow \cite[section 4.2.2.1]{dlebook1} and recent work \cite[Theorem 3.1]{unib}
	which deal with non-degenerate cases on planar domains to work with strong solutions and then by approximations  to establish the following

	\bcoro{2Ndomaindeg}  Let $W$ be a {\em weak} solution to \mref{degsys}. Suppose the spectral-gap condition $\nu_*>1-2/N$ and $\llg(W)\sim|W|^k$ for some $k\ge0$.  Assume also the growth condition
	$$|F(W)|\le C\min\{(|W|^{\frac{k}{2}+2}+1),(|W|^{k+1}+1)\}.$$
	Then, $\|W\|_{BMO(\Og_{2R})}$ is small if $R$ is sufficiently small. 
	
	In particular,  there is a    {\em weak} solution to \mref{degsys} $W$ which is H\"older continuous.

	\ecoro
	
	\bproof  We can replace $dt$ in the proof of 2d case (see \cite[Theorem 3.1]{unib}) by the difference quotient operator $\dg_h$  to deal with weak solutions so that $\ma(W) DW\in L^2(\Og)$. We have
	$$|\ma(W)DW|^2=\myprod{\ma^T(W)\ma(W)DW,DW}\ge \llg^2(W)|DW|^2$$
	so that $\llg(W) DW\in L^2(\Og)$. If $\ma(W)$ is nondegenerate then $DW\in L^2(\Og)$ so that if $N=2$ then $\|W\|_{BMO(\Og_{2R})}$ is small if $R$ is sufficiently small. 
	
	For the degenerate cases, as $\llg(W) DW\in L^2(\Og)$ and therefore $U_i=\int_0^{W_i}\llg(s)ds$ is BMO. Use BMO Johnson lemma 1 we have $U_i^\ag$ is BMO if $0<\ag<1$. Now, as $\llg(W)\sim |W|^k$ with $k>0$ then we take $\ag=(k+1)^{-1}$ to have that $W$ is BMO with $\|W\|_{BMO(\Og_{2R})}$ is small if $R$ is sufficiently small. 
	
	The proof of \reftheo{regdegthm} also applies here for weak solutions in th degenerate/singular case.
	
	Finally, we can consider the strong solutions $W_\eg$'s of the regularized system \mref{degsysregularized}, the above argument then shows that 
	the estimate for $\|W_\eg\|_{BMO(\Og_{2R})}$ is small so that, via a finite covering of $\Og$, $\sup_{t\in(0,t_0)}\|DW_\eg\|_{L^{2p}(\Og)}$ is bounded  independent of $\eg$ for any finite $t_0>0$. There is a subsequence of $\{W_\eg\}$ that converges to a {\em weak} solution $W$ to \mref{degsys} $W$ which is H\"older continuous.
	The proof is complete.
	\eproof
	
	We should note \refrem{strongvsweakGN} concerning the H\"older continuity  of weak solutions. The uniqueness of weak solutions is another issue, if $f$ is Lipschitz in $W$ and $\ma(W)$ is mononicity in some sense then we have uniqueness, see \cite{DB,Vasquez} and \reftheo{uniquethm} below.
	
	The argument in the proof of \refcoro{2Ndomaindeg} cannot be extended to the case $N\ge3$. However, we see that the proof of \cite[Theorem 4.6]{unib} applies once the BMO estimates are established for $N=2$. The corollaries follow \cite[Theorem 4.6]{unib} then apply so that we can state

	\bcoro{thindomainregdeg} Let  $W$ be a {\em weak} solution to \mref{degsys}. Suppose the spectral-gap condition $\nu_*>1-2/N$ and $\llg(W)\sim|W|^k$ for some $k\ge1$.  Assume that
	
	\beqno{BMODxN} \dspl{\int_0^{R}} \iidx{B_R}{|D_{x_N}W(x,x_N,t)|^2}dx_N\le \eg R^{N-2}\eeq
	for all $B_R\subset 	\Og\cap\RR^{N-1}\times\{x_N\}$. Then, $\|W\|_{BMO(\Og_{2R})}$ is small if $R$ is sufficiently small. In particular, there is a    {\em weak} solution to \mref{degsys} $W$ which is H\"older continuous.

	\ecoro

	The condition \mref{BMODxN} is much easier to checked than those in \cite{dlebook,dlebook1} which are directly  involved with estimating $\|DW\|_{L^{2p}}$ ($2p>N$) which is very hard even for strong solutions. Here, we need only consider $|D_{x_N}W|^2$. The arguments lead to \cite[Corollary 4.10]{unib} does not require that $\ma(W)$ is non-degenerate so that we can also consider 'thin' domains in $\RR^N$ for any $N=3$ and have similar existence results as in \cite{unib}.
	
	\section{Uniqueness of weak solutions} \label{unisec}\eqnoset  First of all we adapt the following concept of weak solutions to the equation
	$$\frac{d}{dt}u=\Div(D\Fg(u))+F(x,t)\mbox{ in $\Og\times(0,T)$},\quad u(0)=u_0 \mbox{ on $\Og$}.$$
	
	We say that $u$ is a weak solution on $Q_T=\Og\times(0,T)$ if 
	\bdes
	\item [i)] $u\in L^1(Q_T)$, $\Fg(u)\in L^1(0,T:W^{1,1}_0(\Og))$;
	\item [ii)] for all test function  $\eta\in C^1(\bar{Q_T})$, $\eta=0$ on the parabolic boundary of $Q_T$ we have ($dz=dxdt$)
	$$\itQ{Q_T}{[D\Fg(u)D\eta-u\eta_t]}=\itQ{Q_T}{F\eta}+\iidx{\Og}{u_0\eta(x,0)}.$$
	
	\edes
	
	If $F$ is independent of $u$ then it is well known that there is a weak solution and it is unique (see \cite[section 5.3]{Vasquez}). We can generalize the uniqueness part of  this result if we also assume
	\bdes \item [W)] $u(\cdot,t), D\Fg(u(\cdot,t))\in L^2(\Og)$ for all $t\in(0,T)$ and we can use $u$ as a test function in ii).
	\edes
	
	One can combine Steklov's average and mollifiers to see that W) is reasonable. We then have the following result.
	
	\btheo{uniquethm} Let $\Fg(u)$ be a function on $\RR$, $\mg$ is a function on $\Og\times(0,T)$. Consider  the equation \beqno{degeqn} \frac{d}{dt}u=\Div(D\Fg(u))+\mg(x,t)u,\quad u(0)=u_0 \mbox{ on $\Og$}.\eeq
	Then a weak solution $u$ satisfying W) of this equation is  unique  if \beqno{unicond}\Fg_u(0)=0,\; \Fg_{uu}(tv)v\ge0 \mbox{ for $t\in(0,1)$ and $v\in\RR$},\; \dspl{\int_0^T}\iidx{\Og}{\mg(x,s)}ds<\infty.\eeq
	
	\etheo
	
	The (local) existence of a weak solution of \mref{degeqn} can be established by several means, among these is the approximation method (e.g. see \refrem{existenceweakrem}).
	
	A typical example of \mref{degeqn} with condition \mref{unicond} is: Let $\Fg(u)=|u|^ku$ be a function, $k\ge1$. We have $D\Fg(u)=(k+1)|u|^{k-1}Du$ and $\Fg_{uu}(tv)v=(k+1)k|t|^{k-2}t|v|^{k}\ge0$ for $t\in(0,1)$. The theorem is applied. Also, we can remove W) if one make use of Steklov's average and mollifiers (but we still have to assume $Du(\cdot,t)\in L^2(\Og)$ for all $t\in(0,T)$?). For simplicity, we assume W) here.
	
	Before presenting the proof of \reftheo{uniquethm} let us recall the following well known Gr\"onwall inequality
	\blemm{gronwallineq} Let $y(t), q(t)$ be functions on $(0,T)$. Assume that $y$ is a.e. differentiable and $q(t)$ is integrable. Assume that $$\frac{ d}{dt}y(t)\le q(t)y(t)+c,\; y(0)=y_0\quad\mbox{a.e. on $(0,T)$}, $$
	where $c$ is a constant. Then $y(t)\le e^{\int_0^t q(s)ds}[y_0+c]$.
	
	\elemm
	
	\bproof For any functions $a,b$ (assuming sufficient integrabilities) we have 
	$$\myprod{D(\Fg(a)-\Fg(b)),D(a-b)}=\int_0^1\myprod{\Fg_u(sa+(1-s)b)D(a-b), D(a-b)} ds,$$
	$$\Fg_u(sa+(1-s)b)=\int_0^1 \Fg_{uu}(t(sa+(1-s)b))(sa+(1-s)b)dt+\Fg_u(0).$$
	So that, as $\Fg_u(0)=0$,
	$$\myprod{D(\Fg(a)-\Fg(b)),D(a-b)}=\int_0^1\int_0^1\myprod{ \Fg_{uu}(t(sa+(1-s)b))(sa+(1-s)b)D(a-b), D(a-b)}dt ds.
	$$
	We see that  $\myprod{D(\Fg(a)-\Fg(b)),D(a-b)}\ge0$,as $\Fg_{uu}(tv)v\ge0$ for $t\in(0,1)$ with $v=(sa+(1-s)b)$.

	Consider two weak solutions $u,v$ to the equation \mref{degeqn}. Subtracting, we see that
	$w=u-v$ satisfies
	$\frac{d}{dt}w=\Div(D(\Fg(u)-\Fg(v)))+\mg(x,t)w$. If $u,v$ satisfy W), we can test this equation with  $w$ (using Steklov's average and mollifiers) to obtain
	$$\frac{d}{dt}\iidx{\Og}{|w|^2}+\iidx{\Og}{\myprod{D(\Fg(u)-\Fg(v)),D(u-v)}}= \iidx{\Og}{\mg(x,t)|w|^2}.$$
	As $\myprod{D(\Fg(u)-\Fg(v)),D(u-v)}\ge0$, we can drop the nonegative term on the left hand side to have
	$$\frac{d}{dt}\iidx{\Og}{|w|^2}\le \iidx{\Og}{\mg(x,t)|w|^2},\; w(0)=0$$
	so that by Gr\"onwall 's inequality \reflemm{gronwallineq} and $y_0=c=0$, as  $\dspl{\int_0^T}\iidx{\Og}{\mg(x,s)}ds<\infty$, we have $w\equiv0$ on $\Og\times(0,T)$. \eproof

	\brem{unisysrem} The above calculation can be extended to the system like \mref{degeqn} when $\Fg:\RR^m\to\RR^m$ with $\Fg_{uu}(tv)v$ is a nonnegative definite matrix for $t\in(0,1)$, $v\in\RR^m$. 
	Obviously, we can have $\mg(x,t)=\mbox{diag}[\mg_1(x,t),\ldots,\mg_m(x,t)]$ with $\dspl{\int_0^T}\iidx{\Og}{\mg_i(x,s)}ds<\infty$.
	
	\erem
	
	\brem{degeqnrem} The systems like \mref{degeqn} are also written as $\frac{d}{dt}u=\Div(\hat{\Fg}(u)Du)+q(x,t)u$, with $\hat{\Fg}(u)$ being a $m\times m$ matrix.
	Then a weak solution of this system is unique  if $\hat{\Fg}(0)=0$, $\hat{\Fg}_{u}(tv)v$ is a nonnegative definite matrix for $t\in(0,1)$ and $v\in\RR^m$. The same conditions apply for $\mg_i$'s and we can use approximation to obtain H\"older continuity for weak solutions satisfying W).
	\erem

	{\bf Diagonalize:} We have the following simple result on the boundedness of solutions to a class of cross diffusion systems. Consider the following parabolic system \beqno{diagsys}W_t=\Div(\ma(W) DW).\eeq
	The elliptic case is similar and easier. We can also include a reaction term $f(W)$ with appropriate assumptions in the above systems.
	
	\btheo{boundedsysthm} Assume that there are matrices $B(W),\ag(W)$ and functions $\llg_i(W)$ such that $\ag$ is a diagonal matrix, $\ag(W)=\mbox{diag}[\llg_1(W),\ldots,\llg_m(W)]$, and for some number $\llg_0$
	
	$$B(W)\ma(W) B^{-1}(W)=\ag(W),\; \llg_i(W)\ge\llg_0>0\quad \forall W\in\RR^m.$$ We set 
	$$\mP(W)=\int_0^1 B(sW)Wds.$$ 
	
	If $W$ is a weak solution to \mref{diagsys} and for some constant $c$ (dropping $W$ for abbreviation)\beqno{Bmatcond}\max_i\llg_i|B||(B^{-1})_W||B^{-1}|\le c\min_i\llg_i|\mP|^{-1}\eeq then $\|\mP(W)(\cdot,t)\|_{L^\infty(\Og)}$ is bounded in terms of $\|\mP(W)(\cdot,s)\|_{L^p(\Og)}$ for any $t>s$ and some $p>c+1$.

	\etheo
	
	\bproof We see that $B=\mP_W$ and $DW=B^{-1}D\mP$ and $BW_t=\mP_t$. 
	Multiplying \mref{diagsys} by $B$ and using the fact that $\ma=B^{-1}\ag B$ and the formula  $DB^{-1}=-B^{-1}DB B^{-1}$, we can see that $\mP_t$  is
	$$B\Div(B^{-1}\ag D\mP)=BB^{-1}\Div(\ag D\mP)+B \myprod{D(B^{-1}),\ag D\mP}=\Div(\ag D\mP)- B\myprod{B^{-1}DB B^{-1},\ag D\mP} .$$
	So that $\mP$ satisfies the diagonal system with quadratic growth in gradients ($D\mP$)
	$$\mP_t= \Div(\ag D\mP)- B\myprod{B^{-1}B_W DW B^{-1},\ag D\mP}= \Div(\ag D\mP)- B\myprod{B^{-1}B_W B^{-1} D\mP B^{-1},\ag D\mP}.$$

	We can use the well-known Moser iteration argument by  testing the system with $|\mP|^n\mP$ (using also truncations of the components of $\mP$) to show that $\mP$ is bounded in terms of its $L^p$ norm, which is assumed to be finite. Indeed, we let  $c_*:=|B||B^{-1}B_W B^{-1}||B^{-1}|=|B||(B^{-1})_W||B^{-1}|$ (a function in $W$) and have
	$$\barrl{\frac{1}{n+2}\sup_{t\in (T,T+1)}\iidx{\Og}{|\mP|^{n+2}(x,t)} +(n+1)\int_T^{T+1}\iidx{\Og}{\min_i\llg_i|\mP|^{n-1}|D\mP|^2} \le }{3cm} &\dspl{\int_T^{T+1}}\iidx{\Og}{c_*\max_i\llg_i|\mP|^n|D\mP|^2}+\frac{1}{n+2}\iidx{\Og}{|\mP|^{n+2}(x,T)}.\earr$$
	If $\max_i\llg_i c_*\le c\min_i\llg_i|\mP|^{-1}$ as in \mref{Bmatcond} and $c<n+1$, then the first term on the right hand side can be absorbed to the second term on the left when $n$ large. Standard Moser iteration argument implies (setting $p=n+2$)
	$$\sup_{t\in(T,T+1)} \|\mP(\cdot,t)\|_{L^\infty(\Og)}
	\le C \|\mP(\cdot,T)\|_{L^p(\Og)},\quad \mbox{for  } p>c+1.$$
	
	Note that we can scale $B$ but the number $c$ is unchanged so that we have the above estimate only for some $p>c+1$. The proof is complete. \eproof
	
	This also implies that $\|W\|_{L^\infty}$ is locally bounded in terms of $\|\mP(W)\|_{L^r(\Og)}$ if $\mP^{-1}$ is continuous. 
	
	The easiest example for \mref{Bmatcond} to hold is that $B$ is a constant matrix so that we can take $c=0$ and $n=0$ (so $p=2$). In this case, there is no constraint on $\llg_i(W)$. 
	
	We also note that $|B||B^{-1}|\ge1$ in general.
	If $\min_i\llg_i\sim \max_i\llg_i$ and $B=[b_{ij}]$ and $|b_{ij}|\sim|W|^l$ for some $l\ne-1$ then $|B|\le c_1|W|^l$, $|B^{-1}|\le c_2|W|^{-l}$ and $|(B^{-1})_W|\le lc_2|W|^{-l-1}$. So that  \mref{Bmatcond} is verified as $|B||(B^{-1})_W||B^{-1}|\le lc_1c_2^2|W|^{-(l+1)}\le \frac{lc_1^2}{l+1}c_2^2|\mP|^{-1}$ because, from the definition of $\mP$, we easily see  that $|\mP|\le c_1/(l+1)|W|^{l+1}$. Note that $c=\frac{l}{l+1}c_1^2c_2^2$ which can be very small if $l\sim 0$.
	
	In addition, if $|b_{ij}|\sim|W|^l$ for some $l\ge0$ then  the above results shows that $\|W\|_{L^\infty}$ is locally bounded in terms of $\|W\|_{L^{p(l+1)}(\Og)}$ for  $p\in(c+1,\infty)$. Thus, if $l\sim0$ then $p(l+1)\sim 1$.
	
	\brem{Bmatcondrem} In fact,   we have $B\myprod{B^{-1}B_W B^{-1} D\mP B^{-1},\ag D\mP}= \myprod{B_W B^{-1} D\mP B^{-1},\ag D\mP}$ (as we write $\myprod{x,y}=x\otimes y$ so that $B\myprod{x,y}=Bx\otimes y=\myprod{Bx,y}$) so that the condition \mref{Bmatcond} can be replaced by $\max_i\llg_i|B_W B^{-1}||B^{-1}|\le c\min_i\llg_i|\mP|^{-1}$. If $|b_{ij}|\sim|W|^l$ for some $l\ne-1$ then the number $c$ will put a constraint on the spectral quotient $\min_i\llg_i/ \max_i\llg_i$.
	\erem

	\bibliographystyle{plain}

\end{document}